\batchmode
\documentstyle[12pt,leqno]{amsart}
\newtheorem{dummy}{}[section]
\newtheorem{definition}[dummy]{Definition}
\newtheorem{theorem}[dummy]{Theorem}

\newtheorem{proposition}[dummy]{Proposition}

\newtheorem{corollary}[dummy]{Corollary}
\newtheorem{remark}[dummy]{Remark}

\newtheorem{question}[dummy]{Question}

\begin{document}
\bibliographystyle{plain}
\title{ $p$-adic cocycles and their regulator maps }
\author{Zacky Choo and Victor  Snaith}
\date{March 2009}

%%%%%%%%%%%%%%%%%%%%%%%%%%%%%%%%%%%%%%%%%%%%%%%%
%
%%%%%%%%%%%%%%%%%%%%%%%%%%%%%%%%%%%%%%%%%%%%%%%%

\begin{abstract}
We derive a power series formula for the $p$-adic regulator on the higher dimensional algebraic K-groups of number fields. This formula is designed to be well suited to computer calculations and to reduction modulo powers of $p$. In addition we describe a series of regulator questions concerning higher dimensional K-theoretic analogues of conjectures of Gross and Serre from (\cite{Ta84} Chapter Six).
\end{abstract}
\maketitle
\section{Introduction}

Let $F$ be a $p$-adic local field. Then a $p$-adic regulator is a homomorphism of the form, for $s \geq 2$,
\[     R_{F} :    K_{2s-1}({\cal O}_{F})  \cong  K_{2s-1}(F)  \longrightarrow  F  .   \]
There are lots of such $p$-adic regulators in the literature. They have a number of uses in arithmetic and
geometry. For example, in \cite{Sn97} and \cite{Sn97b} one of us used the $p$-adic regulators together with the
higher dimensional algebraic K-theory local fundamental classes to construct analogues of the classical Chinburg
invariant of the Galois module structure of K-groups. In (\cite{Sn97} \S5) a description is given of five
$p$-adic regulators. These are (a) the cyclic homology regulator (using the results of (\cite{Wag} pp. 244-245;
\cite{Laz}; \cite{LoQu} Theorem 6.2; \cite{LoQu2})  (b)  the dilogarithmic regulator \cite{Col2}  (c) the
\'{e}tale regulator (\cite{DF};  \cite{Sou};  \cite{Sou2} (d) the syntomic regulator (\cite{FM}; \cite{Gr};
\cite{Gr2}; \cite{KTQD})  (e)  the topological cyclic homology regulator (\cite{BeMa}; \cite{HeMa}).  Recently
some of these constructions have been re-examined. In \cite{HK06} Huber and Kings use an idea originally due to
Wagoner \cite{Wag} which substitutes the Lazard isomorphism \cite{Laz} between continuous group cohomology and
Lie algebra cohomology for the van Est isomorphism in the $p$-adic analogue of the construction of the Borel
regulator \cite{BG}. Using the Bloch-Kato exponential \cite{HK06} shows that their $p$-adic regulator coincides
with the \'{e}tale regulator of \cite{Sou2}. A $p$-adic regulator due to Karoubi was overlooked in (\cite{Sn97}
\S5) which, like (a) and (e), also uses cyclic homology and relative K-theory (\cite{Kar83}; \cite{Kar87}). In
\cite{Ham2} the construction of an explicit $p$-adic regulator is sketched which coincides with that of
\cite{Kar83}.

In \cite{CMSS08} we used R.H. Fox's free differential calculus to design an algorithm (implemented in C) to construct explicit homology cycles for the general linear group whose Borel regulators were calculated by power series algorithm (implemented in MAPLE) designed from the explicit formula given in \cite{Ham}. In the course of
working on \cite{CMSS08} and \cite{ZC08} we noticed that our power series also converged $p$-adically, giving rise to an elementary, complete account of the regulator of \cite{Ham2} which culminates in $R_{F}$ of Corollary \ref{4.2}. Independently this construction was introduced in \cite{GT08} and used to show that the regulators of \cite{Ham2} and \cite{HK06} coincide up to a non-zero rational factor.

Our motivation for developing the details of the $p$-adic regulator was similar to our motivation for
\cite{CMSS08}, namely that the power series makes possible an algorithm for evaluating the $p$-adic valuation of the regulators on homology classes in the general linear group of number rings such as those given by the algorithm of \cite{CMSS08}. In all other respects it should be clear to the reader that our approach has nothing to add to the more sophisticated methods of \cite{HK06}, \cite{Kar83} and \cite{GT08}.

\section{Functions on $p$-adic power series}
\begin{definition}
\label{2.0}
\begin{em}

Let $F$ be a $p$-adic local field and let ${\cal O}_{F}$ denote its valuation ring. Let $N$ be a positive integer and let $M_{N}{\cal O}_{F}$ denote the ring of $N \times N$ matrices with entries in ${\cal O}_{F}$ topologised with the $p$-adic topology. Fix a positive odd integer $2s-1$ with $s \geq 2$. Let $E(dx_{0}, \ldots , dx_{2s-1})$
denote the ${\cal O}_{F}$-exterior algebra on symbols $dx_{0}, \ldots , dx_{2s-1}$ so that
$dx_{i} \wedge dx_{j} =  - dx_{j} \wedge dx_{i}$ when $i \not= j$ and $dx_{i} \wedge dx_{i} =0$. Consider the ${\cal O}_{F}$-algebra
\[   \hat{{\cal A} } =   M_{N}{\cal O}_{F}[[x_{0}, x_{1}, \ldots , x_{2s-1}]] \otimes_{{\cal O}_{F}} E(dx_{0}, \ldots , dx_{2s-1})   \]
and set
\[    {\cal A} =   \hat{{\cal A} }/  \simeq ,    \]
the quotient of $ \hat{{\cal A} }$ by the ideal generated by
$1 - \sum_{i=0}^{2s-1} \ x_{i}$ and $\sum_{i=0}^{2s-1} \  dx_{i}$. Setting $| \underline{a} | = a_{0} + a_{1} + \ldots  \ldots + a_{2s-1} $, let $f \in {\cal A}$ have the form
\[   f = \sum_{  \underline{a} = (a_{0}, \ldots , a_{2s-1}) } \ \sum_{u=0}^{2s-1} \
\     f(\underline{a} , u )  p^{e| \underline{a} | }x_{0}^{a_{0}} \cdots x_{2s-1}^{a_{2s-1}} dx_{0} \wedge \cdots \wedge \hat{dx_{u}} \wedge \cdots  dx_{2s-1}   \]
with each $a_{j}  $ an integer greater than or equal to zero and $ f(\underline{a} ,  u )  \in  M_{N}{\cal O}_{F}$.

Define $\Phi_{2s-1}(f)$ by the formula
\[   \Phi_{2s-1}(f)  = \sum_{  \underline{a} = (a_{0}, \ldots , a_{2s-1}) } \ \sum_{u=0}^{2s-1} \  \   (-1)^{u}   {\rm Trace }f(\underline{a} ,  u)  p^{e| \underline{a} | }  \frac{ a_{0}! \cdot  a_{1}! \cdots    a_{2s-1}! }{(| \underline{a} |  + 2s-1)!}  .   \]

Hence $\Phi_{2s-1}(f) $ term-by-term substitutes
\[  (-1)^{u} f(\underline{a} ,  u)  p^{e| \underline{a} | }  \frac{ a_{0}! \cdot  a_{1}! \cdots    a_{2s-1}! }{(| \underline{a} |  + 2s-1)!}   \]
for
\[  f(\underline{a} , u )  p^{e| \underline{a} | }x_{0}^{a_{0}} \cdots x_{2s-1}^{a_{2s-1}} dx_{0} \wedge \cdots \wedge \hat{dx_{u}} \wedge \cdots  dx_{2s-1}  \]
and then takes the trace of the resulting matrix in $M_{N}F$.

By \S\ref{A1} and \S\ref{A2}, $\Phi_{2s-1}(f) $ is well-defined for $f \in {\cal A}$ not merely for $f \in \hat{{\cal A}}$, since it is the term-by-term integral
\[    \Phi_{2s-1}(f)  =  {\rm Trace}  \int_{\Delta^{2s-1}} \  f   .  \]
\end{em}
\end{definition}
\begin{proposition}{$_{}$}
\label{2.1}
\begin{em}

The series $\Phi_{2s-1}(f)$ of Definition \ref{2.0} converges $p$-adically in $F$ for all $e \geq 1$ if $p$ is odd and for all $e \geq 2$ if $p=2$.
\end{em}
\end{proposition}
\vspace{2pt}

Proof:
\vspace{2pt}

If $\nu_{p}$ is the $p$-adic valuation on the rational numbers and $[x]$ denotes the integer part of $x$, then
\[  \nu_{p}( l!)  =   \sum_{i=1}^{\infty} [ \frac{l}{p^{i}}]  =  \frac{ l - \alpha(l)}{p-1} . \]
Here,  if $l =  \sum_{j \geq 0} \  b_{j}p^{j} $ with each $b_{j}$ an integer in the range $0 \leq b_{j} \leq p-1$, we set $\alpha( l) = \sum_{j \geq 0} \  b_{j} $. Therefore
\[  \begin{array}{l}
 \nu_{p}(p^{e| \underline{a} | }  \frac{ a_{0}! \cdot  a_{1}! \cdots  a_{2s-1}! }{(| \underline{a} |  + 2s-1)!}   )   \\
 \\
= e | \underline{a} | + \nu_{p}(a_{0}! \cdot  a_{1}! \cdots  a_{2s-1}! ) -
  \frac{ | \underline{a} |  + 2s-1 -  \alpha( | \underline{a} |  + 2s-1)}{p-1}  \\
  \\
  \geq   (e -  \frac{1}{p-1})   | \underline{a} | -
  \frac{  2s-1}{p-1}
 \end{array}  \]
 which tends to infinity as the monomial $x_{0}^{a_{0}} \cdots x_{2s-1}^{a_{2s-1}}$ tends to zero in the power series topology (that is, as $  | \underline{a} |$ tends to infinity).  $\Box$

\section{A $p$-adic cocycle}
\begin{definition}
\label{3.0}
\begin{em}

Let $X_{0}, X_{1}, \ldots , X_{2s-1}$ be matrices lying in $M_{N}{\cal O}_{F}$ with $s \geq 2$. Denote the $2s$-tuple
$(X_{0}, X_{1}, \ldots , X_{2s-1})$ by $\underline{X}$. If ${\cal A}$ is the algebra introduced in Definition \ref{2.0}, let
\[ \nu(\underline{X}) = 1 + p^{e}  \sum_{i=0}^{2s-1} \  X_{i}x_{i}  \in  {\cal A} .  \]
Hence $\nu(\underline{X})$ is invertible in ${\cal A}$ with
\[   \nu(\underline{X})^{-1} =   1 +  \sum_{i \geq 1} \   (-1)^{i}   p^{e \cdot i}  B^{i}    \]
where $B =   \sum_{i=0}^{2s-1} \  X_{i}x_{i} $. The derivative $d\nu(\underline{X}) =  dB =    \sum_{i=0}^{2s-1} \  X_{i}dx_{i}$ also lies in $ {\cal A}$ and so
\[   \nu(\underline{X})^{-1} d\nu(\underline{X}) \in  {\cal A}.   \]
Furthermore $(\nu(\underline{X})^{-1} d\nu(\underline{X}))^{2s-1}$ is homogeneous of weight $2s-1$ in the differentials $dx_{i}$ so that we have
\[       \Phi_{2s-1}(  (\nu(\underline{X})^{-1} d\nu(\underline{X}))^{2s-1}) \in F   . \]

Denote by $G_{N, e}F$ the closed subgroup of $GL_{N}{\cal O}_{F}$ consisting of matrices which are congruent to the identity modulo $p^{e}$. With the $p$-adic topology on  $G_{N, e}F$ the map
\[  \tilde{\Phi}_{2s-1} :   (1 + p^{e}X_{0}, 1 + p^{e}X_{1}, \ldots , 1 + p^{e}X_{2s-1})   \mapsto
 \Phi_{2s-1}(  (\nu(\underline{X})^{-1} d\nu(\underline{X}))^{2s-1})  \]
lies in ${\rm Map}_{cts}(   (G_{N, e}F)^{2s} , F )  $,  the $p$-adically continuous functions from the $2s$-fold cartesian product of $G_{N, e}F$ to $F$.
\end{em}
\end{definition}
\begin{theorem}{$_{}$}
\label{3.1}
\begin{em}

(i)  \   If  $Y_{1}, Y_{2} \in  G_{N,e}F$ then
\[   \begin{array}{l}
 \tilde{\Phi}_{2s-1}(Y_{1}(1 + p^{e}X_{0})Y_{2},  \ldots , Y_{1}(1 + p^{e}X_{2s-1})Y_{2})  \\
 \\
 =
  \tilde{\Phi}_{2s-1}(1 + p^{e}X_{0},  \ldots , 1 + p^{e}X_{2s-1})  .
  \end{array}   \]

  Similarly, if $Y \in GL_{N}{\cal O}_{F}$,
  \[   \begin{array}{l}
 \tilde{\Phi}_{2s-1}(Y(1 + p^{e}X_{0})Y^{-1},  \ldots , Y(1 + p^{e}X_{2s-1})Y^{-1})  \\
 \\
 =
  \tilde{\Phi}_{2s-1}(1 + p^{e}X_{0},  \ldots , 1 + p^{e}X_{2s-1})  .
  \end{array}   \]

  (ii)  \    If $F/E$ is a Galois extension and $\sigma \in {\rm Gal}(F/E)$ then
\[   \sigma(   \tilde{\Phi}_{2s-1}(1 + p^{e}X_{0},  \ldots , 1 + p^{e}X_{2s-1}) ) =   \tilde{\Phi}_{2s-1}(1 + p^{e} \sigma X_{0},  \ldots , 1 + p^{e} \sigma X_{2s-1}) . \]

(iii)  \   The function $  \tilde{\Phi}_{2s-1}$ is a $(2s-1)$-dimensional $p$-adically continuous cocycle on $G_{N,e}F$ with values in the trivial $G_{N,e}F$-module $F$.

\end{em}
\end{theorem}
\vspace{2pt}

{\bf Proof}
\vspace{2pt}

For part (i) we have, in Definition \ref{3.0},
\[   \begin{array}{l}
   \nu(\underline{X})^{-1} d\nu(\underline{X})  \\
   \\
   =   ( \sum_{i=0}^{2s-1} \  (1 + p^{e}X_{i})x_{i} ))^{-1} d ( \sum_{i=0}^{2s-1} \  (1 + p^{e}X_{i})x_{i} )
\end{array}  \]
while
\[   \begin{array}{l}
   ( \sum_{i=0}^{2s-1} \  Y_{1}(1 + p^{e}X_{i})Y_{2}x_{i} )^{-1} d ( \sum_{i=0}^{2s-1} \  Y_{1}(1 + p^{e}X_{i})Y_{2}x_{i} )  \\
   \\
   =   Y_{2}^{-1} ( \sum_{i=0}^{2s-1} \  (1 + p^{e}X_{i})x_{i} ))^{-1} d ( \sum_{i=0}^{2s-1} \  (1 + p^{e}X_{i})x_{i} )Y_{2}   \\
   \\
 =   Y_{2}^{-1}  \nu(\underline{X})^{-1} d\nu(\underline{X}) Y_{2} .
\end{array}  \]
Hence the first part of (i) follows from the integral formula of Definition \ref{2.0}
\[  \begin{array}{ll}
\Phi_{2s-1}( (  \nu(\underline{X})^{-1} d\nu(\underline{X}))^{2s-1})    \\
\\
  =  {\rm Trace}  \int_{\Delta^{2s-1}} \
( \nu(\underline{X})^{-1} d\nu(\underline{X}))^{2s-1}   \\
\\
 =   {\rm Trace}  \int_{\Delta^{2s-1}} \
Y_{2}^{-1}  ( \nu(\underline{X})^{-1} d\nu(\underline{X}))^{2s-1} Y_{2}  .
\end{array} \]
The proof of the second part of (i) is similar.

Part (ii) is immediate from Definition \ref{3.0}.

For Part (iii) we must first prove that $\Phi_{2s-1}$ lies in ${\rm Map}_{cts, G_{N,e}F}( (G_{N,e}F)^{2s} , F)$ where $G_{N,e}F$ acts trivially on $F$ and by diagonal  left multiplication on $(G_{N,e}F)^{2s}$. This follows from the first part of (i).

We use the elementary form Stokes' Theorem from \S\ref{A3} (cf. \cite{Ham}) to prove the cocycle condition.
Given a $(2s+1)$-tuple of matrices $(X_{0}, X_{1}, X_{2}, \ldots , X_{2s})$ in $M_{N}{\cal O}_{F}$ form
$\nu = 1 +   \sum_{i=0}^{2s} \  p^{e} X_{i}x_{i}$ where the $x_{i}$ are the barycentric coordinates in $\Delta^{2s}$ and corresponding to the $i$-th face for each $0 \leq i \leq 2s$ we set
\[ \begin{array}{l}
\nu_{i} =  1+  p^{e} X_{0}x_{0}+  p^{e} X_{1}x_{1}+ p^{e} X_{2}x_{2}+  \ldots  \hat{p^{e} X_{i}x_{i}} +  \ldots + p^{e} X_{2s}x_{2s}  ,
\end{array}  \]
deleting the $i$-th term from $\nu$.
Then the cocycle condition is the vanishing of the expression
\[    \begin{array}{c}
 \sum_{i=0}^{2s} \    (-1)^{i} {\rm Trace} \int_{(x_{i}=0) \bigcap \Delta^{2s}}  \  (\nu_{i}^{-1} d\nu_{i})^{2s-1} .
\end{array}  \]

Select a monomial $(2s-1)$-form from within $(\nu^{-1} d\nu)^{2s-1}$, say
\[   \omega  = x_{0}^{a_{0}} x_{1}^{a_{1}}  x_{2}^{a_{2}}  \ldots  x_{2s-1}^{a_{2s-1}} x_{2s}^{a_{2s}}   dx_{0} \wedge \ldots \wedge \hat{dx_{u}}  \wedge \ldots \wedge \hat{dx_{v}} \wedge \ldots dx_{2s} \]
with $0 \leq u < v \leq 2s$. By \S \ref{A3}
 \[ \begin{array}{l}
 \sum_{i=0}^{2s} \  (-1)^{i} \    \int_{(x_{i}=0) \bigcap \Delta^{2s}}  \   \omega   =
  \int_{\Delta^{2s}}  \  d \omega
 \end{array} \]
 so that, since the sums of these identities converge $p$-adically, by Proposition \ref{2.1},
 \[    \begin{array}{c}
 \sum_{i=0}^{2s} \    (-1)^{i} {\rm Trace} \int_{(x_{i}=0) \bigcap \Delta^{2s}}  \  (\nu_{i}^{-1} d\nu_{i})^{2s-1} =   {\rm Trace}  \int_{\Delta^{2s}}  \   d  (\nu^{-1} d\nu)^{2s-1}  .
\end{array}  \]
However $0 = d(\nu^{-1} \nu) =  d(\nu^{-1}) \nu + \nu^{-1} d \nu$ so that
\[  d  (\nu^{-1} d\nu)^{2s-1}  =  - (2s-1)   (\nu^{-1} d\nu)^{2s}  \]
which implies that
\[   {\rm Trace}  \int_{\Delta^{2s}}  \   d  (\nu^{-1} d\nu)^{2s-1} = 0 \]
because the integrand changes sign under the $2s$-cycle permutation of the $\nu^{-1} d\nu$'s and a cyclic permutation of a product of matrices preserves the trace. $\Box$

\section{The $p$-adic regulator}
\begin{definition}
\label{4.0}
\begin{em}

As in \S \ref{2.0}, let $F$ be a $p$-adic local field and let ${\cal O}_{F}$ denote its valuation ring. From the localisation sequence for algebraic K-theory \cite{Qu73} and the vanishing of even K-groups of finite fields \cite{Qu72} we have an isomorphism
\[    K_{2s-1}({\cal O}_{F})   \cong K_{2s-1}(F)   \]
for all $s \geq 2$. Let ${\rm Hur}$ denote the Hurewicz homomorphism to the integral homology of the infinite general linear group, with the discrete topology,
\[  {\rm Hur} :    K_{2s-1}({\cal O}_{F})  \longrightarrow  H_{2s-1}(GL{\cal O}_{F} ; {\mathbb Z} ). \]
When $N$ is large the inclusion induces an isomorphism
\[   H_{2s-1}(GL_{N}{\cal O}_{F} ; {\mathbb Z} )   \stackrel{\cong}{\longrightarrow}  H_{2s-1}(GL{\cal O}_{F} ; {\mathbb Z} )  .  \]
To be precise this is true for $N \geq {\rm max}(4s-1 ,  2s-1 + sr({\cal O}_{F}))$ where $sr({\cal O}_{F})$
is Bass's stable rank of ${\cal O}_{F}$ \cite{NS89}.

Choosing $N$ large we define
\[    R_{N, F} :  H_{2s-1}( GL_{N}{\cal O}_{F} ; {\mathbb Z} )   \longrightarrow  F   \]
to be equal to the composition of the transfer map
\[    H_{2s-1}( GL_{N}{\cal O}_{F} ; {\mathbb Z} )
 \stackrel{Tr}{\longrightarrow}       H_{2s-1}( G_{N,e}F ; {\mathbb Z} )   \]
with the homomorphism
\[   \frac{1}{ [ GL_{N}{\cal O}_{F} : G_{N,e}F  ]}   \cdot  <  [ \tilde{\Phi}_{2s-1}] , - >  :
  H_{2s-1}( G_{N,e}F ; {\mathbb Z} )   \longrightarrow  F   \]
  given by pairing a discrete homology class with the continuous cohomology class of Theorem 3.1(iii)
and dividing by the index of $ G_{N,e}F $ in $GL_{N}{\cal O}_{F}$. Explicitly, if $d, \epsilon$ are the residue degree and ramification index of $F / {\mathbb Q}_{p}$
\[  [ GL_{N}{\cal O}_{F} : G_{N,e}F  ] =  |GL_{N}{\mathbb F}_{p^{d}}| p^{N^{2} d (e \epsilon - 1)}  . \]
\end{em}
\end{definition}
\begin{proposition}{$_{}$}
\label{4.1}
\begin{em}

For large $N$ the homomorphism $ R_{N, F} $ is equal to the homomorphism
\[    H_{2s-1}(GL_{N}{\cal O}_{F} ; {\mathbb Z} )   \stackrel{i_{*}}{\longrightarrow}  H_{2s-1}(GL_{N+1}{\cal O}_{F} ; {\mathbb Z} )   \stackrel{R_{N+1 , F}}{\longrightarrow}  F      \]
where $i_{*}$ is induced by the inclusion map.
\end{em}
\end{proposition}
\vspace{2pt}

{\bf Proof}
\vspace{2pt}

Dualising the Double Coset Formula of (\cite{SnEBI} p.19) we have the homology version for $J, H \subseteq G$, subgroups of finite index.
\[  \begin{array}{l}
  H_{2s-1}(J ; {\mathbb Z} )   \stackrel{i_{*}}{\longrightarrow}  H_{2s-1}(G ; {\mathbb Z} )   \stackrel{Tr}{\longrightarrow}       H_{2s-1}( H ; {\mathbb Z} )  \\
  \\
  = \sum_{z \in  J \backslash   G  /  H } \
H_{2s-1}( J ; {\mathbb Z} )    \stackrel{Tr}{\longrightarrow}   H_{2s-1}( J \bigcap zHz^{-1} ; {\mathbb Z})  \\
\\
\hspace{30pt}  \stackrel{(z^{-1} - z)_{*}}{\longrightarrow}
  H_{2s-1}( z^{-1}Jz \bigcap H ; {\mathbb Z} )   \stackrel{i_{*}}{ \longrightarrow }
    H_{2s-1}(  H ; {\mathbb Z} ) .
    \end{array}         \]
We wish to apply this to the case in which $J = GL_{N}{\cal O}_{F}$, $G =  GL_{N+1}{\cal O}_{F}$
and $H = G_{N+1,e}F \lhd G$. In this case $zHz^{-1} = H $ and so
\[  \begin{array}{l}
  H_{2s-1}(J ; {\mathbb Z} )   \stackrel{i_{*}}{\longrightarrow}  H_{2s-1}(G ; {\mathbb Z} )   \stackrel{Tr}{\longrightarrow}       H_{2s-1}( H ; {\mathbb Z} )  \\
  \\
  = \sum_{z \in  J \backslash   G  /  H } \
H_{2s-1}( J ; {\mathbb Z} )    \stackrel{Tr}{\longrightarrow}   H_{2s-1}( J \bigcap H ; {\mathbb Z})  \\
\\
\hspace{30pt}   \stackrel{i_{*}}{ \longrightarrow }
    H_{2s-1}(  H ; {\mathbb Z} )
  \stackrel{(z^{-1} - z)_{*}}{\longrightarrow}   H_{2s-1}(  H ; {\mathbb Z} )   .
    \end{array}         \]

From the second part of (i)
\[      Res_{G_{N,e}F}^{G_{N+1,e}F}( (z^{-1} - z)^{*}[ \Phi_{2s-1}])  \in  H_{cts}^{2s-1}( G_{N,e}F; {\mathbb Z})   \]
is equal to $[\Phi_{2s-1}]$ for $G_{N,e}F$. Therefore
\[  \begin{array}{l}
[GL_{N+1}{\cal O}_{F} :  G_{N+1,e}F ] R_{N+1 , F} \cdot i_{*}   \\
\\
=   [GL_{N+1}{\cal O}_{F} :  GL_{N}{\cal O}_{F} \cdot  G_{N+1,e}F ] [GL_{N}{\cal O}_{F} ; G_{N,e}F ]  R_{N,F}  \\
\\
=    [GL_{N+1}{\cal O}_{F} :  GL_{N}{\cal O}_{F} \cdot  G_{N+1,e}F ] [GL_{N}{\cal O}_{F} G_{N+1,e}F ; G_{N+1,e}F ]  R_{N,F}
\end{array}      \]
so that $R_{N,F} = R_{N+1 , F} \cdot i_{*}$.  $\Box$
\begin{corollary}{$_{}$}
\label{4.2}
\begin{em}

For large $N$ the homomorphism
\[   H_{2s-1}(GL{\cal O}_{F}; {\mathbb Z})   \stackrel{i_{*}^{-1}}{\longrightarrow }
 H_{2s-1}(GL_{N}{\cal O}_{F}; {\mathbb Z})   \stackrel{R_{N,F}}{\longrightarrow }   F   \]
 is independent of $N$.
\end{em}
\end{corollary}
\begin{definition}
\label{4.3}
\begin{em}

Define a homomorphism
\[    \hat{R}_{F}  :  H_{2s-1}(GL{\cal O}_{F}; {\mathbb Z})  \longrightarrow    F   \]
by the formula
\[     \hat{R}_{F}   =     \frac{(-1)^{s} (s-1)!}{(2s-2)! (2s-1)!}  R_{N,F} i_{*}^{-1}  ,  \]
in the notation of Corollary \ref{4.2}, where $N$ is a large positive integer.
\end{em}
\end{definition}
\begin{theorem}{$_{}$}
\label{4.4}
\begin{em}

In the notation of Definitions \ref{4.0} and \ref{4.3} the composition
\[    R_{F} :  K_{2s-1}(F)    \cong     K_{2s-1}({\cal O}_{F})  \stackrel{ {\rm Hur}}{\longrightarrow}
H_{2s-1}(GL{\cal O}_{F} ; {\mathbb Z} )   \stackrel{ \hat{R}_{F}  }{\longrightarrow}   F  \]
is equal to the $p$-adic regulator homomorphism defined in \cite{Ham2} (and hence also with that of \cite{HK06}).
\end{em}
\end{theorem}
\vspace{2pt}

{\bf Proof}
\vspace{2pt}

First we should point out that \cite{Ham2} gives an explicit formula for a $p$-adic regulator only in the case when $F = {\mathbb Q}_{p}$. However the sketched proof showing that this construction is well-defined and coincides with Karoubi's cyclic homology $p$-adic regulator applies equally well for general $F$. The regulator of \cite{Ham2} is defined by composing the Hurewicz homomorphism with the homomorphism, for large $N$,
\[   R :   H_{2s-1}(GL_{N}F ; {\mathbb Z} )   \longrightarrow   F  \]
which is induced by sending a $2s$-tuple of matrices $(Y_{0}, \ldots , Y_{2s-1})$ in the bar resolution for $GL_{N}F$ to the integral
\[        \frac{(-1)^{s} (s-1)!}{(2s-2)! (2s-1)!} {\rm Trace} \int_{\Delta^{2s-1}}  \  (\nu^{-1} d \nu)^{2s-1}  \]
where $\nu = \sum_{i=0}^{2s-1} x_{i}Y_{i}$  where the $x_{i}$'s are the barycentric coordinates. The verification that this integral converges $p$-adically for a general $2s$-tuple is quite delicate and is carried out in the Appendix to \cite{GT08}.

On the other hand, the construction which we have given uses the same integral, but only in the situation where each $Y_{i} $ lies in $G_{N,e}F$ in which case we saw in \S2 and \S3 that it is very easy to show  $p$-adic convergence.

Let $j :  G_{N,e}F  \longrightarrow  GL_{N}{\cal O}_{F}$ denote the inclusion. The above discussion shows that
\[   [ GL_{N}{\cal O}_{F} : G_{N,e}F  ]  R_{N,F} =  R \cdot j_{*} \cdot Tr  :
  H_{2s-1}( GL_{N}{\cal O}_{F} ; {\mathbb Z} )   \longrightarrow  F    \]
and the result follows since $ j_{*} \cdot Tr  =  [ GL_{N}{\cal O}_{F} : G_{N,e}F  ] $. $\Box$
\begin{remark}
\label{4.5}
\begin{em}

Using an explicit $p$-adically analytic cocycle it is shown in \cite{GT08} that the Karoubi-Hamida $p$-adic regulator, which equals $R_{F}$ by Theorem \ref{4.4}, coincides up to a non-zero rational factor with the Wagoner-Huber-Kings $p$-adic regulator of \cite{Wag} and \cite{HK06}.
\end{em}
\end{remark}

 \section{K-theoretic Analogies of (\cite{Ta84} Chapter Six)}

 The construction of the homomorphism $R_{F} $ of Theorem \ref{4.4} makes sense when $s=1$ providing that we restrict attention to $K_{1}({\cal O}_{F}) \cong {\cal O}_{F}^{*}$ rather than
 $K_{1}(F) \cong F^{*}$. Taking $N=1$ then, as in \S\ref{4.0},
\[  [ GL_{1}{\cal O}_{F} : G_{1,e}F  ] =  [ {\cal O}_{F}^{*} :  1 + p^{e}{\cal O}_{F} ] =  (p^{d} - 1) p^{ d (e \epsilon - 1)}  \]
where $d, \epsilon$ are the residue degree and ramification index of $F / {\mathbb Q}_{p}$ respectively.
When $s=1$ in the constructions of Definitions \ref{2.0} and \ref{3.0} yield
\[   \tilde{\Phi}_{1}( 1 + p^{e} X_{0} , 1 + p^{e} X_{1})  =  {\rm log}_{p}( \frac{1 + p^{e} X_{1}}{1 + p^{e} X_{0}}) \]
for $X_{0}, X_{1} \in {\cal O}_{F}$ where ${\rm log}_{p}( 1 + p^{e}z) = \sum_{i=0}^{\infty} \  (-1)^{i} \frac{(p^{e}z)^{i+1}}{i+1}$ is the usual $p$-adic logarithm series. Therefore when $s=1$ one finds that
\[ R_{F} : {\cal O}_{F}^{*} \longrightarrow  F \]
is given by $R_{F}(x) = - {\rm log}_{p}(x)$ where, as usual, $ {\rm log}_{p}$ is the unique homomorphic extension of the $p$-adic logarithm series to ${\cal O}_{F}^{*}$. Similarly, the Borel regulator on the units in a ring of algebraic integers of a number field is essential equal to the Archimdean logarithm \cite{Ham}. Therefore, just as there are analogs of the Stark conjectures involving the Borel regulators
on higher dimensional algebraic K-groups of number fields (\cite{Sn04}, \cite{Sn05}, \cite{Sn06}), so there are higher dimensional algebraic K-theoretic analogs of the $p$-adic versions  of Stark's conjecture, due to B. Gross and J-P. Serre, described in (\cite{Ta84} Chapter Six).

Accordingly, in this section we shall examine possible analogs of the material in  (\cite{Ta84} Chapter Six) involving higher K-groups. We shall begin with a simple reciprocity law which hardly features logarithms (i.e. regulators) at all!
\begin{dummy}{$p$-adic absolute values}
\label{7.1}
%\label{1.1}
\begin{em}

Let $k$ be a number field, $v$ a place of $k$ and $p$ a place of ${\mathbb Q}$ (if $p = \infty$ then ${\mathbb Q}_{p} = {\mathbb R}$). Let $x \in k^{*}$. If $v$ is a finite place
\[    |x|_{v} =  (Nv)^{-v(x)}   \in {\mathbb Q}  . \]
For $v$ Archimedean we set
\[   (Nv)^{v(x)}   =  \left\{
\begin{array}{ll}
1  &  {\rm if}  \  v \  {\rm is \ complex},  \\
\\
{\rm sign}(\sigma(x)) & {\rm if} \   v   \ {\rm is \  real, \ induced \ by} \   \sigma : k \rightarrow  {\mathbb R}.
\end{array}  \right.  \]
Recall from (\cite{Ta84} p.7) that $Nv$  and the normalised absolute values are defined in terms of the $v$-adic completion $k_{v}$ by:  $|-|_{v} $ is the usual absolute value when $k_{v} \cong {\mathbb R}$, $|x|_{v} = x \overline{x}$ when $k_{v} \cong {\mathbb C}$ and when $v$ is discrete with uniformiser $\pi$ then $Nv = | {\cal O}_{k}/(\pi) |$ and $| \pi |_{v} = (Nv)^{-1}$.

Therefore we have
\[    {\rm Norm}_{k/{\mathbb Q}}(x)   =   \prod_{v} \  (Nv)^{v(x)}  . \]
\end{em}
\end{dummy}
\begin{definition}
\label{7.2}
%\label{1.2}
\begin{em}

Define the ${\mathbb Q}_{p}$-valued absolute value of $x$ to be the element $|x|_{v,p} \in {\mathbb Q}_{p}^{*}$ given by
\[  |x|_{v,p} = \left\{  \begin{array}{ll}
   {\rm Norm}_{k_{v}/{\mathbb Q}_{p}}(x) \cdot   (Nv)^{-v(x)}  & {\rm if} \  v | p ,  \\
   \\
    (Nv)^{-v(x)}  &   {\rm if} \  v  \not|  p
\end{array}  \right.   \]
Hence $ |x|_{v,p} = \pm 1$ when $v$ is Archimedean.
\end{em}
\end{definition}
%\begin{remark}
\label{7.3}
%\label{1.3}
%\begin{em}

(a)  \   For $v = \infty$ we have  $ |x|_{v, \infty} =  |x|_{v} \in {\mathbb R}$ and if $p$ is finite  $|x|_{v,p}   \in {\mathbb Z}_{p}^{*}$ for all $x \in k^{*}$.

(b)  \   For all $x \in k^{*}$
\[   \prod_{v} \   |x|_{v,p}    = 1   \in   {\mathbb Q}_{p}  \]
and
\[     \sum_{v} \   {\rm log}_{p}(|x|_{v,p})    = 0   \in   {\mathbb Q}_{p}.  \]
For finite primes $v$ we have
\[   \prod_{v \  {\rm finite}} \   |x|_{v,p}    =  \pm 1   \in   {\mathbb Q}_{p}  \]
and
\[     \sum_{v  \  {\rm finite}  } \   {\rm log}_{p}(|x|_{v,p})    = 0   \in   {\mathbb Q}_{p}.  \]
(c)  \   When $p$ is a prime number $|-|_{v,p}$ coincides with the composition
\[   k^{*}  \stackrel{i}{\longrightarrow} k_{v}^{*}  \stackrel{recip}{\longrightarrow}  {\rm Gal}(k_{v}^{ab}/k_{v})
\stackrel{\chi}{\longrightarrow}   {\mathbb Z}_{p}^{*}   \]
where $i$ is the inclusion, $recip$ is the reciprocity map and $\chi$ is the inverse of the cyclotomic character giving the action of ${\rm Gal}(k_{v}^{ab}/k_{v})  $ on the $p$-primary roots of unity.
\begin{dummy}{An analogue of Remark \ref{7.3}(b)}
\label{7.4}
\begin{em}

Let $k$ be a number field, let $v$ a finite place of $k$ and $p$ be a rational prime. Let $x \in  K_{2s-1}(k)$ for some $s \geq 2$. Recall that $K_{2s-1}(k)$ is a finitely generated abelian group isomorphic to $K_{2s-1}({\cal O}_{k})$ where ${\cal O}_{k}$ is the ring of algebraic integers of $k$.
There is a canonical higher dimensional local fundamental class (\cite{Sn97}, \cite{Sn97b}) associated to a Galois extension of the form $L/k_{v}$. This is an element of ${\rm Ext}_{{\mathbb Z}[{\rm Gal}(L/k_{v})]}^{2}( {\rm Tors}K_{2s-2}(L),  K_{2s-1}(L))$ which is represented by a $2$-extension
\[    K_{2s-1}(L)  \longrightarrow  A  \longrightarrow  B  \longrightarrow   {\rm Tors}K_{2s-2}(L)    \]
with $A, B$ cohomologically trivial. Therefore we have a canonical reciprocity isomorphism
between Tate cohomology groups
\[    \hat{H}^{0}( {\rm Gal}(L/k_{v}) ; K_{2s-1}(L) )
\stackrel{\cong}{\longrightarrow}  \hat{H}^{-2}( {\rm Gal}(L/k_{v}) ;  {\rm Tors}K_{2s-2}(L)) \]
which may be identified (\cite{SnF} Definition 1.1.2,  p.3) with
\[    \frac{K_{2s-1}(k_{v})}{ {\rm Norm} K_{2s-1}(L) }
\stackrel{\cong}{\longrightarrow}    H_{1}( {\rm Gal}(L/k_{v}) ;  {\rm Tors}K_{2s-2}(L)) .  \]

Recall from (\cite{HS71} Chapter VI, \S4) that if $G$ is a finite group and $M$ is a ${\mathbb Z}[G]$-module we have an isomorphism
\[  {\mathbb Z} \otimes_{ {\mathbb Z}[G]}  IG \cong  IG/(IG)^{2}  \cong  G^{ab}  \]
given by $1 \otimes_{G} (g-1)   \mapsto    (g-1) \ (modulo   \   IG^{2})   \mapsto  g[G,G] \in  G^{ab}$.
More generally, we have a short exact sequence of left $G$-modules with the diagonal action
\[   0  \longrightarrow  IG \otimes M   \longrightarrow  {\mathbb Z}[G] \otimes M    \longrightarrow   M   \longrightarrow  0 \]
where the right-hand map is $g \otimes m \mapsto $ and the resulting long exact homology sequence looks like
\[    0 = H_{1}(G ;   {\mathbb Z}[G] \otimes M  )  \longrightarrow   H_{1}(G ;    M  ) \longrightarrow
{\mathbb Z} \otimes_{G}   IG \otimes M  \longrightarrow  {\mathbb Z} \otimes_{G}    {\mathbb Z}[G] \otimes M  \cong M   \]
where the right-hand map is $1 \otimes_{G}  (g-1) \otimes m \mapsto  (g^{-1}-1)m$. Hence
\[   H_{1}(G ;    M  )   \cong {\rm Ker}(  {\mathbb Z} \otimes_{G}   IG \otimes M  \longrightarrow  M). \]

Given a homomorphism $\chi :G \longrightarrow  {\mathbb Z}_{p}^{*}$ there is an induced homomorphism
\[   \chi' :   {\mathbb Z} \otimes_{G}   IG \otimes M  \longrightarrow    {\mathbb Z}_{p}^{*} \otimes  \frac{M}{IG \cdot M } \]
given by $1 \otimes (g-1) \otimes m  \mapsto  \chi(g) \otimes  (m \  {\rm modulo} \  IG \cdot M) $. This is well-defined because
\[  \begin{array}{l}
\chi'( 1 \otimes g' (g-1) \otimes g' m ) \\
\\
=     \chi(g' g) \otimes  (g' m \  {\rm modulo} \  IG \cdot M)
-    \chi(g') \otimes  (g' m \  {\rm modulo} \  IG \cdot M)   \\
\\
=    \chi(g') \otimes  (g' m \  {\rm modulo} \  IG \cdot M) +  \chi(g) \otimes  (g' m \  {\rm modulo} \  IG \cdot M) \\
\\
\hspace{70pt}  -  \chi(g') \otimes  (g' m \  {\rm modulo} \  IG \cdot M)  \\
\\
=   \chi(g) \otimes  (g' m \  {\rm modulo} \  IG \cdot M) \\
\\
=  \chi(g) \otimes   (m \  {\rm modulo} \  IG \cdot M) .
\end{array}   \]
Setting ${\rm Gal}(L/k_{v}) $ and $M =  {\rm Tors}K_{2s-2}(L)$ the inverse of the cyclotomic character
induces a canonical homomorphism
\[    H_{1}( {\rm Gal}(L/k_{v}) ;  {\rm Tors}K_{2s-2}(L)) \longrightarrow    {\mathbb Z}_{p}^{*} \otimes
 \frac{ {\rm Tors}K_{2s-2}(L)}{   I{\rm Gal}(L/k_{v})   \cdot  {\rm Tors}K_{2s-2}(L) } .  \]
In addition, in the course of the proof of naturality for the higher dimensional fundamental classes (\cite{Sn97}, \cite{Sn97b}) it is shown that there is a natural isomorphism of the form
 \[  \frac{ {\rm Tors}K_{2s-2}(L)}{ I{\rm Gal}(L/k_{v})   \cdot  {\rm Tors}K_{2s-2}(L) } \cong    {\rm Tors}K_{2s-2}(k_{v})  . \]
 Composing these homomorphism when $s \geq 2$ yields a homomorphism
 \[   \chi(s)_{v,p}  : K_{2s-1}(k)  \longrightarrow   K_{2s-1}(k_{v})   \longrightarrow   {\mathbb Z}_{p}^{*} \otimes    {\rm Tors}K_{2s-2}(k_{v})   \]
 which, when $s=1$, is equal to
 \[    |-|_{v,p}  : K_{1}(k) = k^{*} \longrightarrow   {\mathbb Z}_{p}^{*} .  \]
 Since there is an isomorphism of $k$-algebras of the form
 \[     k \otimes_{{\mathbb Q}}  {\mathbb Q}_{p}  \longrightarrow   \prod_{v|p}  \   k_{v}  \]
 we obtain a homomorphism
\[   \{    \chi(s)_{v,p}  \}_{v|p}  :    K_{2s-1}(k)  \longrightarrow    {\mathbb Z}_{p}^{*} \otimes    {\rm Tors}K_{2s-2}(  k \otimes_{{\mathbb Q}}  {\mathbb Q}_{p} ) .  \]
The Galois behaviour of this homomorphism and the reciprocity of \S\ref{7.3}(b) suggests the following question :
\begin{question}
\label{7.5}
\begin{em}

For $s \geq 2$ does the image of $ \{    \chi(s)_{v,p}  \}_{v|p}$ lie in the image of the homomorphism
\[   K_{2s-2}( k)   \longrightarrow     {\mathbb Z}_{p}^{*} \otimes    {\rm Tors}K_{2s-2}(  k \otimes_{{\mathbb Q}}  {\mathbb Q}_{p} ) \]
induced by the inclusion of $k$? Perhaps this is true with $ K_{2s-2}( k) $ replaced by $ K_{2s-2}( {\cal O}_{k}) $, the K-group of the ring of integers of $k$?
\end{em}
\end{question}
\end{em}
\end{dummy}

\begin{dummy}{Higher dimensional Stark conjectures}
\label{7.6}
\begin{em}

For the reader's convenience let us recall the analogue for higher dimensional algebraic K-theory of the classical Stark conjecture of (\cite{Ta84} Chapter One, \S5). This conjecture was posed by one of us in \cite{Sn04}, \cite{Sn05} and \cite{Sn06}, unaware that B.H. Gross \cite{Gross} had already asked this question about the Stark conjecture decades earlier in the preprint, which eventually appeared as \cite{Gross05}.

Let $K/k$ be a Galois extension of number fields.
Let $\Sigma(K)$ denote the set of embeddings of $K$ into the  complex numbers.
For $r = -1, -2, -3, \ldots$ set
\[ Y_{r}(K) = \prod_{\Sigma(K)} \ (2 \pi i )^{-r} {\mathbb Z}  =  {\rm Map}( \Sigma(K) , (2 \pi i )^{-r} {\mathbb Z} )  \]
endowed with the ${\rm Gal}(\mathbb{C}/{\mathbb{R}})$-action diagonally on $\Sigma(K)$ and on $(2 \pi i )^{-r}
$and $Y_{r}(K)^{+}$ denotes the subgroup fixed by complex conjugation. Therefore
\[ rank_{{\bf Z}}( Y_{r}(K)^{+} ) = \left\{
\begin{array}{ll}
r_{2} & {\rm if } \   r  \  {\rm is \ odd}, \\
\\
r_{1} + r_{2} & {\rm if } \   r   \  {\rm is \ even}   .
 \end{array} \right. \]
where $| \Sigma(K)|  = r_{1} + 2r_{2}$ and $r_{1}$ is the number of real embeddings of $K$.
Denote by ${\cal O}_{K}$ the integers of $K$.
For any negative integer $r $ we have the Borel regulator (\cite{BG}, \cite{HW1})
\[ R_{K}^{r} : K_{1-2r}( {\cal O}_{K} ) \otimes {\mathbb R}
\stackrel{\cong}{\longrightarrow}
  Y_{r}(K)^{+} \otimes {\mathbb R}  \]
which is an ${\mathbb R}[{\rm Gal}(K/k)]$-isomorphism. Choose a
${\mathbb Q}[{\rm Gal}(K/k)]$-isomorphism of the form
\[ f_{r,K } : K_{1-2r}( {\cal O}_{K} ) \otimes {\mathbb Q}
\stackrel{\cong}{\longrightarrow}
  Y_{r}(K)^{+} \otimes {\mathbb Q}  \]
so that
\[ R_{K}^{r} \cdot (f_{r, K})^{-1} :   Y_{r}(K)^{+} \otimes {\mathbb R}
\stackrel{\cong}{\longrightarrow}
  Y_{r}(K)^{+} \otimes {\mathbb R}  \]
is an ${\mathbb R}[{\rm Gal}(K/k)]$-isomorphism. Then we form the Stark regulator defined,
for each representation $V$ of ${\rm Gal}(K/k)$, by
\[  R(V, f_{r,K}) = det( (R_{K}^{r} \cdot f_{r, K}^{-1})_{*} \in Aut_{{\mathbb
C}}(Hom_{{\rm Gal}(K/k)]}( V^{\vee} ,  Y_{r}(K)^{+}
\otimes {\mathbb C}) ))  , \]
where $V^{\vee}$ is the contragredient representation of $V$.

Let $S$ be a finite set of primes of $k$ which includes all the Archimdedean primes and all the finite primes which ramify in $K/k$. Let $L_{k,S}^{*}(r , V) $ denote the leading term of the Taylor
expansion of the Artin L-function associated to $S$ and $V$ at $s = r$.
We define a function
$ {\cal R}_{f_{r,K}}$ given on a finite-dimensional complex representation $V$ by
\[ {\cal R}_{f_{r,K}}(V) = \frac{R(V,f_{r,K})}{L_{k,S}^{*}(r , V) }.  \]

Then the higher-dimensional analogue of the Stark conjecture of \cite{Ta84}
asserts that, if $\Omega_{{\mathbb Q}}$ denotes the absolute Galois group of the rationals,
\[  {\cal R}_{f_{r,K}} \in Hom_{\Omega_{{\mathbb Q}}}( R( {\rm Gal}(K/k)] ) ,
\overline{{\mathbb Q}}^{*})  \subseteq Hom( R( {\rm Gal}(K/k)] ) ,{\mathbb C}^{*})  \]
and the truth of this conjecture is independent of the choice of $f_{r,K}$.

The calculations of Beilinson (\cite{Bei}; see also \cite{BGreit} \S4.2, \cite{HW1} and \cite{Neuk}) show that the
higher-dimensional analogue of the Stark conjecture is true when $K/k$ is a subextension of
any abelian extension of the rationals (see \cite{Sn06} Theorem 7.6 (proof)).
\end{em}
\end{dummy}
\begin{dummy}{$p$-adic L-functions}
\label{7.7}
\begin{em}

Let $\overline{k}$ denote an algebraic closure of $k$. If $p$ is a prime, let
\[   \omega : {\rm Gal}( \overline{k} /k) \longrightarrow   \mu({\mathbb Q}_{p})   \]
denote the Teichm\"{u}ller character (\cite{Ta84} p.130) associated with the Galois action on the $p$-power roots of unity and taking values in the $p$-adic roots of unity $\mu({\mathbb Q}_{p}) $.

Let ${\mathbb C}_{p}$ denote the $p$-adic completion of an algebraic closure of ${\mathbb Q}_{p}$ (\cite{Ta84} p.129).
Let $V$ be a continuous, finite-dimensional representation of $ {\rm Gal}( \overline{k} /k) $ over ${\mathbb C}_{p}$. Suppose that $\{ 1 , \tau \}$ is the decomposition group of a place of $\overline{k}$ whose restriction to $k$ is real. Following (\cite{Ta84} p.130) we shall call such an element $\tau$ a ``conjugation''. Then $V$ is totally even if any such $\tau$ acts trivially on $V$ and is totally odd if any such $\tau$ acts as minus the identity.

Now let $V$ be a finite-dimensional ${\mathbb C}_{p}$-representation of ${\rm Gal}(K/k)$ where $K/k$ is a finite
Galois extension. Therefore, for all integers $n$, $V \otimes \omega^{n}$ is also a ${\mathbb
C}_{p}$-representation of ${\rm Gal}( \overline{k} /k)$ which factors through the Galois group of a finite
extension of $k$.

Now let $\alpha : {\mathbb C}_{p} \stackrel{\cong}{\longrightarrow} {\mathbb C}$ denote an isomorphism of fields so that, if ${\rm dim}_{{\mathbb C}_{p} }(V) = t$, we may form the complex representation
 \[   \alpha(V \otimes \omega^{n}) :  {\rm Gal}(K/k)  \longrightarrow  GL_{t}{\mathbb C}   \]
by choosing a matrix representation of $V$ and applying $\alpha$ to the matrix entries.

Let $S$ be a finite set of primes of $k$ which includes all the Archimedean  primes and all the finite primes which divide $p$. Then the $p$-adic L-function is the unique meromorphic function (\cite{Ta84} p.131)
\[   L_{p, S}(-, V)  :  {\mathbb Z}_{p} \longrightarrow  {\mathbb C}_{p}  \]
which satisfies the interpolation formula
\[    \alpha( L_{p, S}(n , V))   =    L_{k, S}( n  ,  \alpha(V \otimes \omega^{n-1}))   \]
for all strictly negative integers $n$ and all field isomorphisms $\alpha : {\mathbb C}_{p} \stackrel{\cong}{\longrightarrow} {\mathbb C}$. The functional equation for the Artin L-function (\cite{Ta84} p.20) together with Euler's functional equation for the $\gamma$-function
\[       \Gamma(s) =   \frac{\Gamma(s+n+1)}{z(z+1) \cdots (z+n)}  \]
shows that $ L_{p, S}(-, V) $ is identically zero unless $k$ is totally real and $V$ is totally even. Constructions of the $p$-adic L-function are given in \cite{Cas-Nog} and \cite{DR80} (see also \cite{Bar} and \cite{Quee}).
\end{em}
\end{dummy}
\begin{dummy}{$p$-adic Higher dimensional Stark conjectures}
\label{7.8}
\begin{em}

The Stark conjecture at $s=0$ features the Dirichlet regulator constructed from the Archimedean logarithm and the $p$-adic analogue at $s=0$, due to Gross (\cite{Ta84} p.132) replaces the logarithm by the $p$-adic logarithm. As remarked at the beginning of this section, the $p$-adic regulator in $K_{1}$ is minus the $p$-adic logarithm and the Borel regulator behaves similarly. Therefore it is natural to formulate similar conjectures on higher dimensional K-groups using their $p$-adic regulator maps, $R_{F}$ of Theorem \ref{4.4}.

Let $K/k$ be a finite Galois extension of number fields with $k$ totally real and $K$ totally imaginary.

Consider $K \otimes_{{\mathbb Q}} F$ where $F/ {\mathbb Q}_{p}$ is an extension of local fields. We have $K = {\mathbb Q}(\beta)$ for some algebraic $\beta$ whose minimal polynomial is $m_{\beta}(x) \in {\mathbb Q}[x]$. Suppose that
$m_{\beta}(x)$ splits in $F$ then we have
\[    K \otimes_{{\mathbb Q}} F  \cong   {\mathbb Q}[x]/ (m_{\beta}(x)) \otimes_{{\mathbb Q}} F  \cong
F[x]/ m_{\beta}(x) \cong    \prod_{i=1}^{{\rm deg}(m_{\beta}(x)) }  \   F  \]
where the last map evaluates polynomials at each of the distinct roots of $m_{\beta}(x)$. Therefore the composition
\[    K  \longrightarrow    K \otimes_{{\mathbb Q}} F   \stackrel{\cong}{\longrightarrow}   \prod_{1}^{[K : {\mathbb Q}]} \ F     \]
is given by $ z \mapsto z \otimes 1 \mapsto  \{  z_{i}  \}  $ where $z_{i}$ is the image of $z$ under the inclusion of $K$ into $F$ corresponding to the $i$-th root of $m_{\beta}(x)$.

Take the case where $F = {\mathbb C}_{p}$ then we have an involutive field automorphism
\[    c_{p} =  \alpha \cdot c \cdot \alpha^{-1} :  {\mathbb C}_{p} \stackrel{\cong}{\longrightarrow}   {\mathbb C}
\stackrel{\cong}{\longrightarrow}   {\mathbb C}_{p}    \]
where $c$ is complex conjugation. This depends upon the choice of $\alpha$. The analogue of the diagonal action of $c$ on $Y_{r}(K) \otimes {\mathbb R}$ is the involution on $ \prod_{1}^{[K : {\mathbb Q}]} \ F $ which sends $F$ in the coordinate corresponding to the $i$-th root $w_{i}$ by $c_{p}$ to the copy of $F$ corresponding to the root $c_{p}(w_{i})$.

Let $\sigma_{i} : K    \longrightarrow {\mathbb C}_{p} $ denote the embedding corresponding to the $i$-th root $w_{i}$. Let $Y_{(p)}(K) = \prod_{\Sigma_{p}(K)}  \  {\mathbb Z}  $ where $\Sigma_{p}(K)$ is the set of embeddings of $K$ into ${\mathbb C}_{p}$. Define a homomorphism
\[   R^{r}_{p, K} : K_{1-2r}(K)  \longrightarrow  Y_{(p)}(K) \otimes_{{\mathbb Z}} {\mathbb C}_{p}   \]
to have the $\sigma_{i}$-th coordinate given by the composition
\[   K_{1-2r}(K)  \stackrel{(\sigma_{i})_{*}}{\longrightarrow}  K_{1-2r}({\mathbb C}_{p})
 \stackrel{ R_{ {\mathbb C}_{p} }}{\longrightarrow}  {\mathbb C}_{p}  .    \]
 By naturality of the $p$-adic regulator this is equal to the composition
\[   K_{1-2r}(K)  \stackrel{(\sigma_{i})_{*}}{\longrightarrow}  K_{1-2r}({\mathbb Q}_{p}(w_{i}))
 \stackrel{ R_{ {\mathbb Q}_{p}(w_{i})}}{\longrightarrow} {\mathbb Q}_{p}(w_{i})  \longrightarrow  {\mathbb C}_{p}  .   \]
\end{em}
\end{dummy}

\begin{proposition}{$_{}$}
\label{7.10}
\begin{em}

In \S\ref{7.8}  $ R_{p,K}^{r}$ is a ${\rm Gal}(K/k)$-homomorphism whose image lies in $(Y_{(p)}(K) \otimes_{{\mathbb Z}} {\mathbb C}_{p} )^{+} $, the $(+1)$-eigenspace of the involution $c_{p}$.
\end{em}
\end{proposition}
\vspace{2pt}

{\bf Proof}
\vspace{2pt}

Suppose that the root $w_{2i}$ is the complex conjugate of $w_{2i-1}$ then, by Galois equivariance of the $p$-adic regulator,  an
element in the image of $R_{p, K}^{r} $ will have the $(2i-1 , 2i)$-th pair of coordinates of the form
$(z_{2i-1} , c_{p}(z_{2i-1}))$. However this pair is sent by the involution to
$(c_{p}(z_{2i-1}) , c_{p}(c_{p}(z_{2i-1})))= (z_{2i-1} , c_{p}(z_{2i-1}))$, as required.    $\Box $

\begin{dummy}{The original $p$-adic Gross conjecture}
\label{7.11}
\begin{em}

Let $p$ be a fixed prime. Let $S$ be a finite set of places of $K$ containing all the Archimedean places and the places dividing $p$. Let ${\cal O}_{K,S}^{*} \cong K_{1}({\cal O}_{K,S})$ denote the $S$-units of $K$. Set $Y = \oplus_{v \in S}  \  {\mathbb Z}$ denote the free abelian group on the elements of $S$ and set
\[    X =  \{   \sum_{v \in S} \  n_{v} \cdot v  \in  Y  \    |   \   \   \sum_{v \in S} \  n_{v} = 0  \}  .   \]
Define
\[     \lambda_{p}  :   K_{1}({\cal O}_{K,S})   \longrightarrow   {\mathbb Q}_{p}  \otimes_{{\mathbb Z}}   X  \]
by $   \lambda_{p}(x) =    \sum_{v \in S} \   log_{p}(  |x|_{v,p})  \cdot  v$ where $ |x|_{v,p} \in {\mathbb Z}_{p}^{*}$ is as in Definition \ref{7.2}. We also denote by $ \lambda_{p}  $ the linear extension
\[     \lambda_{p}  :    {\mathbb C}_{p}  \otimes_{{\mathbb Z}}  K_{1}({\cal O}_{k,S})   \longrightarrow   {\mathbb C}_{p}  \otimes_{{\mathbb Z}}   X  .  \]
Let $f ; X \longrightarrow    {\mathbb C}_{p}  \otimes_{{\mathbb Z}}  K_{1}({\cal O}_{K,S}) $ be a ${\rm Gal}(K/k)$-homomorphism. For all $\alpha : {\mathbb C}_{p} \stackrel{\cong}{\longrightarrow} {\mathbb C}$ we set
\[   \alpha(f) =  (\alpha \otimes 1) \cdot f :  X \longrightarrow   {\mathbb C}   \otimes_{{\mathbb Z}}  K_{1}({\cal O}_{K,S})   .  \]

Define
\[    A_{p , 0}(f , V) =   \frac{{\rm det}( 1 \otimes (\lambda_{p} \cdot f) , (V \otimes {\mathbb C}_{p}  \otimes_{{\mathbb Z}}   X)^{{\rm Gal}(\overline{k}/k)})}{  L_{p, S}^{*}(0 , V \otimes \omega ) } \]
where $ L_{p, S}^{*}(n , V \otimes \omega )$ is the leading term of the Taylor series for the $p$-adic L-function of $k$ at $s=n$ and $V$ is a totally odd representation of ${\rm Gal}(K/k)$.

The $p$-adic Gross conjecture asserts that
\[   \alpha(A_{p , 0}(f , V)) =   {\cal R}_{\alpha(f)}(\alpha(V))    \]
where $ {\cal R}_{\alpha(f)}(\alpha(V)) $ is the analogue at $s=0$ of the non-zero complex number defined in \S\ref{7.6} using the leading term of the $L$-function at $s=  -1, -2, -3, \ldots$.
\end{em}
\end{dummy}
\begin{dummy}{The higher dimensional $p$-adic conjecture}
\label{7.12}
\begin{em}

As is usual in the conjecturing business, we slavishly follow the earlier conjecture making a systematic change. In this case the change is to replace $ log_{p}(  |x|_{v,p}) $ by $  R^{r}_{p, K} $.

Following Gross, for $r=  -1, -2, -3, \ldots$, we now assume we are given a ${\rm Gal}(K/k)$-homomorphism
\[   f :  Y_{(p)}(K)  \longrightarrow   K_{1-2r}(K) \otimes_{{\mathbb Z}}  {\mathbb C}_{p} \]
which me extend linearly to give a ${\mathbb C}_{p}[{\rm Gal}(K/k)]$-homomorphism
\[   f :  Y_{(p)}(K)    \otimes_{{\mathbb Z}}  {\mathbb C}_{p}    \longrightarrow   K_{1-2r}(K) \otimes_{{\mathbb Z}}  {\mathbb C}_{p}    \]
and then form the composition
\[   R^{r}_{p, K}   \cdot  f^{+} :   (Y_{(p)}(K)    \otimes_{{\mathbb Z}}  {\mathbb C}_{p})^{+}
  \longrightarrow   K_{1-2r}(K) \otimes_{{\mathbb Z}}  {\mathbb C}_{p}
  \longrightarrow   (Y_{(p)}(K)    \otimes_{{\mathbb Z}}  {\mathbb C}_{p})^{+} .  \]

  Then the higher dimensional analogue of \S \ref{7.11} would assert that
  \[    A_{p , r}(f , V) =   \frac{{\rm det}( 1 \otimes (   R^{r}_{p, K}  \cdot f^{+}) , (V \otimes ( {\mathbb C}_{p}  \otimes_{{\mathbb Z}}    Y_{(p)}(K))^{+})^{{\rm Gal}(\overline{k}/k)})}{  L_{p, S}^{*}(r , V \otimes \omega^{1-r} ) }   \]
satisfies
\[    \alpha(A_{p , r}(f , V)) =   {\cal R}_{\alpha(f)}(\alpha(V))   \]
when $V \otimes \omega^{1-r}$ is totally even and $  {\cal R}_{\alpha(f)}(\alpha(V)) $ is as in \S\ref{7.6} at
$s = r$.
\end{em}
\end{dummy}
\begin{proposition}{$_{}$}
\label{7.13}
\begin{em}

Suppose that $K/k$ is a subextension of an abelian extension over ${\mathbb Q}$. Then if the
higher dimensional $p$-adic conjecture of \S\ref{7.12} is true for one $\alpha$ then it is true for all $\alpha$'s.
\end{em}
\end{proposition}
\vspace{2pt}

{\bf Proof}
\vspace{2pt}

The higher dimensional Stark conjecture of \S\ref{7.6} is true for cyclotomic extensions and their subextensions (\cite{Sn06} \S 3.1; see also \cite{Sn04} and \cite{Sn05}).
Therefore the result follows by the argument of (\cite{Ta84}, Chapter Six Theorem 5.2).   $\Box$
\begin{remark}
\label{7.14}
\begin{em}

The truth of the conjecture of \S\ref{7.12} implies that the $p$-adic regulator
\[   R^{r}_{p, K} : K_{1-2r}(K)  \otimes_{{\mathbb Z}} {\mathbb C}_{p}  \longrightarrow  (Y_{(p)}(K) \otimes_{{\mathbb Z}} {\mathbb C}_{p})^{+}  \]
is an isomorphism because the determinant of $ R^{r}_{p, K}   \cdot  f^{+}$ is non-zero and the
$ {\mathbb C}_{p}$-dimensions of  domain and range are both equal to $r_{2}$.

In a subsequent paper we shall verify that $ R^{r}_{p, K}$ is indeed an isomorphism.
\end{em}
\end{remark}

\section{Appendix: Elementary integration}
\begin{dummy}{Scalar integration}
\label{A1}
\begin{em}

Let $\Delta^{n} =  \{  (x_{0}, \ldots , x_{n}) \in {\mathbb R}^{n+1} \ | \  0 \leq x_{i} \ {\rm and } \  \sum_{i=0}^{n} x_{i}=1 \} $ denote the usual $n$-simplex. Consider the iterated integral
\[    \int_{x_{n} = 0 }^{\rho(n)}   \  \int_{x_{n-1} = 0 }^{\rho(n-1)}  \
\cdots    \int_{x_{0} = 0 }^{\rho(0)}   \   f  dx_{0} \ldots  \hat{dx_{i}} \ldots   dx_{n}   \]
where each $a_{j} \geq 0$, the integral corresponding to $x_{i}$ is omitted,
\[  f = x_{0}^{a_{0}}  x_{1}^{a_{1}}  \cdots (  1 -  \ldots x_{i-1} - \hat{x_{i}}  - x_{i+1} \ldots )^{a_{i}} \cdots  x_{n}^{a_{n}}  ,  \]
$\rho(n) = 1$,  $\rho(n-1) = 1-x_{n} $ and in general,
\[   \rho(j) = 1 - x_{j+1} - x_{j+2} - \ldots - \hat{x_{i}} - \ldots  - x_{n}.   \]
We have
\[    \begin{array}{l}
   \int_{x_{0}=0}^{\rho(0) }  \    x_{0}^{a_{0}}  (\rho(0) -  x_{0})^{a_{i}}   dx_{0}  \\
   \\
   =   \left\{    \begin{array}{ll}
   \rho(0)  &  {\rm if} \   a_{0}=0=a_{i}  ,  \\
   \\
   \frac{\rho(0)^{a_{0}+1}}{(a_{0}+1) }   \  & {\rm if } \  a_{0}>0,  a_{i}=0 , \\
   \\
     \frac{\rho(0)^{a_{i}+1}}{(a_{i}+1) }   \   &  {\rm if } \  a_{i}>0,  a_{0}=0 , \\
   \\
        \int_{x_{0}=0}^{\rho(0) }  \  a_{0} x_{0}^{a_{0}-1}  \frac{(\rho(0) - x_{0})^{a_{i} + 1} }{(a_{i} + 1)} dx_{0} &   {\rm if}   \    a_{0}, a_{i} >0
   \end{array}  \right.   \\
   \\
   =        \frac{  a_{0}! a_{i}!  \rho(0)^{a_{0}+a_{i}}}{(a_{0}+ a_{i}+1)!}  \\
   \\
   =     \frac{  a_{0}! a_{i}! ( \rho(1)-x_{1})^{a_{0}+a_{i}}}{(a_{0}+ a_{i}+1)!}  ,
   \end{array}   \]
   integrating by parts. By induction we find that, for $a_{i} \geq 0$,
\[    \begin{array}{l}
   \int_{x_{n} = 0 }^{\rho(n)}   \  \int_{x_{n-1} = 0 }^{\rho(n-1)}  \
\cdots    \int_{x_{0} = 0 }^{\rho(0)}   \   x_{0}^{a_{0}}  x_{1}^{a_{1}} \ldots x_{n}^{a_{n}}dx_{0} \ldots  \hat{dx_{i}} \ldots   dx_{n}   \\
\\
=  \frac{ a_{0}! \cdot  a_{1}! \cdots  \cdots  a_{n-1}! \cdot  a_{n}! }{(a_{0} + a_{1} + \ldots  \ldots + a_{n} + n)!}  .
   \end{array}   \]

\end{em}
\end{dummy}
\begin{dummy}{Integration of differential forms}
\label{A2}
\begin{em}

The $n$-simplex
\[ \Delta^{n} = \{ (x_{0}, \ldots , x_{n}) \  |  \   0 \leq x_{i}, \ \sum_{i} \ x_{i}=1  \}  \subset {\mathbb R}^{n+1}  \]
is an orientable, $n$-dimensional manifold with boundary. For the purposes of integration of a differentiable $n$-form
\[  f(x_{0}, \ldots , x_{n}) dx_{0} \wedge  \ldots \wedge \hat{dx_{v}} \wedge  \ldots \wedge dx_{n}   \]
on $ \Delta^{n}$ we use
\[  \begin{array}{l}
 dx_{0} \wedge \ldots \wedge \hat{dx_{v}} \wedge  \ldots \wedge dx_{n}   \\
 \\
 =
( - \sum_{j=1}^{n} \ dx_{j})  \wedge dx_{1} \wedge \ldots \wedge \hat{dx_{v}} \wedge  \ldots \wedge dx_{n}   \\
\\
=  -   dx_{v}   \wedge dx_{1}  \wedge  \ldots \wedge \hat{dx_{v}} \wedge  \ldots \wedge dx_{n}   \\
\\
=    (-1)^{v}   dx_{1}   \wedge  \ldots \wedge dx_{v} \wedge  \ldots \wedge dx_{n}
 \end{array}  \]
 to rewrite
 \[ \int_{\Delta^{n}} \   f  dx_{0} \wedge  \ldots \wedge \hat{dx_{v}} \wedge  \ldots \wedge dx_{n}  =
 (-1)^{v}   \int_{\Delta^{n}} \   f  dx_{1}  \wedge  \ldots \wedge dx_{v} \wedge  \ldots \wedge dx_{n} . \]
The embedding $\phi_{n}$ of ${\mathbb R}^{n}$ into ${\mathbb R}^{n+1}$ given by
\[    \phi_{n}(y_{1}, \ldots , y_{n}) = (1 - \sum_{i} y_{i} , y_{1}, \ldots , y_{n})  \]
maps
\[ \underline{ \Delta}^{n}  =    \{ (y_{1}, \ldots , y_{n}) \  |  \   0 \leq y_{i}, \ \sum_{i} \ y_{i}  \leq 1  \}  \subset {\mathbb R}^{n}  \]
diffeomorphically onto $\Delta^{n}$. We define (see \cite{BT} p.31)
\[ \begin{array}{l}
\int_{\Delta^{n}} \   f  dx_{0} \wedge  \ldots \wedge \hat{dx_{v}} \wedge  \ldots \wedge dx_{n}  \\
\\
=  (-1)^{v}     \int_{y_{n} = 0 }^{\rho(n)}   \  \int_{y_{n-1} = 0 }^{\rho(n-1)}  \
\cdots    \int_{y_{1} = 0 }^{\rho(1)}   \   f(\phi_{n}^{-1}(y_{1}, \ldots , y_{n}) )  dy_{1}  \ldots   dy_{n}
\end{array} \]
where $\rho(n) = 1$,  $\rho(n-1) = 1-y_{n} $ and in general,
\[   \rho(j) = 1 - y_{j+1} - y_{j+2} - \ldots   - y_{n}.   \]
In other words the integral on $\Delta^{n}$ is transformed to an integral on $\underline{ \Delta}^{n} $ with respect to the standard volume form $dy_{1} \wedge y_{2} \wedge  \ldots \wedge dy_{n}$ on ${\mathbb R}^{n}$. In particular we find, if each integer $a_{j}$ is greater than or equal to zero, by \S\ref{A1}
\[  \begin{array}{l}
 \int_{\Delta^{n}} \     x_{0}^{a_{0}}  x_{1}^{a_{1}}  \cdots x_{i}^{a_{i}} \cdots  x_{n}^{a_{n}}  dx_{0} \wedge  \ldots \wedge \hat{dx_{v}} \wedge  \ldots \wedge dx_{n}  \\
\\
=  (-1)^{v}   \frac{ a_{0}! \cdot  a_{1}! \cdots  \cdots  a_{n-1}! \cdot  a_{n}! }{(a_{0} + a_{1} + \ldots  \ldots + a_{n} + n)!}   .
\end{array} \]
\end{em}
\end{dummy}
\begin{dummy}{Stokes' Theorem for monomial differential forms}
\label{A3}
\begin{em}

Suppose that $n=2s$ and that we have a monomial $(2s-1)$-form on $\Delta^{2s}$
\[   \omega  = x_{0}^{a_{0}} x_{1}^{a_{1}}  x_{2}^{a_{2}}  \ldots  x_{2s-1}^{a_{2s-1}} x_{2s}^{a_{2s}}   dx_{0} \wedge \ldots \wedge \hat{dx_{u}}  \wedge \ldots \wedge \hat{dx_{v}} \wedge \ldots \wedge dx_{2s} \]
with $0 \leq u < v \leq 2s$ and each $a_{j}$ an integer greater than or equal to zero. Hence the differential $d \omega$ is given by the expression
\[   \begin{array}{l}
  (\sum_{j=0}^{2s}  \   a_{j} x_{0}^{a_{0}}  \ldots   x_{j}^{a_{j}-1}  \ldots  x_{2s-1}^{a_{2s-1}} x_{2s}^{a_{2s}}  dx_{j} ) \wedge \ldots \wedge \hat{dx_{u}}  \wedge \ldots \wedge \hat{dx_{v}} \wedge \ldots \wedge dx_{2s}   \\
\\
=    a_{u} x_{0}^{a_{0}}  \ldots   x_{u}^{a_{u}-1}  \ldots  x_{2s-1}^{a_{2s-1}} x_{2s}^{a_{2s}}  dx_{u} \wedge  dx_{0} \wedge \ldots \wedge \hat{dx_{u}}  \wedge \ldots \wedge \hat{dx_{v}} \wedge \ldots \wedge dx_{2s} \\
\\
\hspace{5pt} +
 a_{v} x_{0}^{a_{0}}  \ldots   x_{v}^{a_{v}-1}  \ldots  x_{2s-1}^{a_{2s-1}} x_{2s}^{a_{2s}}  dx_{v} \wedge  dx_{0} \wedge \ldots \wedge \hat{dx_{u}}  \wedge \ldots \wedge \hat{dx_{v}} \wedge \ldots \wedge dx_{2s}  \\
 \\
 =   (-1)^{u}  a_{u} x_{0}^{a_{0}}  \ldots   x_{u}^{a_{u}-1}  \ldots  x_{2s-1}^{a_{2s-1}} x_{2s}^{a_{2s}}     dx_{0} \wedge \ldots \wedge dx_{u}  \wedge \ldots \wedge \hat{dx_{v}} \wedge \ldots \wedge dx_{2s} \\
\\
\hspace{5pt} +
(-1)^{v+1} a_{v} x_{0}^{a_{0}}  \ldots   x_{v}^{a_{v}-1}  \ldots  x_{2s-1}^{a_{2s-1}} x_{2s}^{a_{2s}}   dx_{0} \wedge \ldots \wedge \hat{dx_{u}}  \wedge \ldots \wedge dx_{v} \wedge \ldots \wedge  dx_{2s}
\end{array}  \]
so that, by \S\ref{A2},
\[   \int_{\Delta^{2s}}  \  d \omega  =   \left\{  \begin{array}{ll}
0 &  {\rm if}  \  a_{u}=0, a_{v}=0 ,  \\
\\
  (-1)^{u+v+1}  \frac{ a_{0}! \cdot  a_{1}! \cdots  \cdots  a_{2s-1}! \cdot  a_{2s}! }{(a_{0} + a_{1} + \ldots  \ldots + a_{2s} + 2s-1)!}  & {\rm if} \   a_{u}=0, a_{v} > 0  \\
  \\
 (-1)^{u+v}     \frac{ a_{0}! \cdot  a_{1}! \cdots  \cdots  a_{2s-1}! \cdot  a_{2s}! }{(a_{0} + a_{1} + \ldots  \ldots + a_{2s} + 2s-1)!}  & {\rm if} \   a_{v}=0, a_{u} > 0  \\
      \\
        0   & {\rm if} \   a_{v} > 0, a_{u} > 0 .
\end{array} \right.    \]

Now consider the restriction of $\omega$ to the $(2s-1)$-simplex $(x_{i}=0) \bigcap \Delta^{2s}$. The integral
\[   \int_{(x_{i}=0) \bigcap \Delta^{2s}}  \   \omega \]
is zero unless $a_{i}=0$ and $i \in \{ u , v \}$. The ordered coordinates for $(x_{i}=0) \bigcap \Delta^{2s}$ are $(x_{0}, x_{1}, \ldots , x_{i-1}, x_{i+1} , \ldots , x_{2s})$ so that, by \S\ref{A2},
\[  \begin{array}{l}
 \int_{(x_{i}=0) \bigcap \Delta^{2s}}  \   \omega \\
 \\
 =   \left\{    \begin{array}{ll}
 (-1)^{v+1}   \frac{ a_{0}! \cdot  a_{1}! \cdots  \cdots  a_{2s-1}! \cdot  a_{2s}! }{(a_{0} + a_{1} + \ldots  \ldots + a_{2s} + 2s-1)!}    & {\rm if}  \   a_{i}=0, \  i = u , \\
 \\
 (-1)^{u}  \frac{ a_{0}! \cdot  a_{1}! \cdots  \cdots  a_{2s-1}! \cdot  a_{2s}! }{(a_{0} + a_{1} + \ldots  \ldots + a_{2s} + 2s-1)!}    & {\rm if}  \   a_{i}=0, \  i = v ,  \\
 \\
 0 & {\rm otherwise} .
 \end{array}  \right.
 \end{array}   \]
 Therefore
 \[ \begin{array}{l}
 \sum_{i=0}^{2s} \  (-1)^{i} \    \int_{(x_{i}=0) \bigcap \Delta^{2s}}  \   \omega   \\
 \\
 =    \left\{  \begin{array}{ll}
((-1)^{u+v+1} + (-1)^{u+v})  \frac{ a_{0}! \cdot  a_{1}! \cdots  \cdots  a_{2s-1}! \cdot  a_{2s}! }{(a_{0} + a_{1} + \ldots  \ldots + a_{2s} + 2s-1)!}  &  {\rm if}  \  a_{u}=0, a_{v}=0 ,  \\
\\
  (-1)^{u+v+1}  \frac{ a_{0}! \cdot  a_{1}! \cdots  \cdots  a_{2s-1}! \cdot  a_{2s}! }{(a_{0} + a_{1} + \ldots  \ldots + a_{2s} + 2s-1)!}  & {\rm if} \   a_{u}=0, a_{v} > 0  \\
  \\
 (-1)^{u+v}     \frac{ a_{0}! \cdot  a_{1}! \cdots  \cdots  a_{2s-1}! \cdot  a_{2s}! }{(a_{0} + a_{1} + \ldots  \ldots + a_{2s} + 2s-1)!}  & {\rm if} \   a_{v}=0, a_{u} > 0  \\
      \\
        0   & {\rm if} \   a_{v} > 0, a_{u} > 0 .
\end{array} \right.  \\
\\
=  \int_{\Delta^{2s}}  \  d \omega .
 \end{array} \]
\end{em}
\end{dummy}
\section{Appendix: Explicit formulae for the transfer}

Suppose that $H$ is a subgroup of $G$ of finite index, $[G : H] = m$ and let $\{ x_{i} \ | \ 1 \leq i \leq m \}$ be a set of right coset representatives for $H \backslash G$. Let $B_{*}G$ denote the bar resolution. Hence $B_{n}G$ is the free abelian group on $G^{n+1}$ for $n \geq 0$ with differential
$d : B_{n}G \longrightarrow B_{n-1}G$ given by
\[    d(g_{0}, \ldots , g_{n}) =  \sum_{j=0}^{n} \ (-1)^{j} \ (g_{0} , \ldots , \hat{g_{j}} , \ldots , g_{n}) \]
and left $G$-module structure given by $g( g_{0}, \ldots , g_{n}) = (gg_{0}, \ldots , gg_{n}) $. Hence
$H_{*}(G ; {\mathbb Z})$ is given by the homology of the chain complex
\[  \ldots  \stackrel{ 1 \otimes d}{\longrightarrow}   {\mathbb Z} \otimes_{{\mathbb Z}[G]} B_{n}G  \stackrel{ 1 \otimes d}{\longrightarrow}     {\mathbb Z} \otimes_{{\mathbb Z}[G]} B_{n-1}G   \stackrel{ 1 \otimes d}{\longrightarrow}   \ldots  .  \]

There is an anti-homomorphism to the symmetric group
\[ \pi : G \longrightarrow  \Sigma_{m}  \]
given by the right action of $G$ on $H \backslash G$
\[    x_{i} g =    h(i,g) x_{\pi(g)(i)}   \]
for a unique $h(i,g)  \in H$. We have
\[    x_{i} g g'=    h(i,g) x_{\pi(g)(i)} g'  =   h(i,g) h( \pi(g)(i) , g')  x_{ \pi(g')(\pi(g)(i))}  \]
so that $\pi(g g') =   \pi(g') \cdot \pi(g) \in  \Sigma_{m}$ and
\[  h(i,g)^{-1} h(i, g g') =   h( \pi(g)(i) , g')  \in H .  \]

Since $ {\mathbb Z} \otimes_{{\mathbb Z}[G]} B_{n}G$ is the free abelian group on $\{ 1 \otimes_{G} (1, g_{1}, \ldots , g_{n}) \  |  \  g_{i}  \in G \}$ we may define homomorphisms
\[  T_{n} :   {\mathbb Z} \otimes_{{\mathbb Z}[G]} B_{n}G  \longrightarrow   {\mathbb Z} \otimes_{{\mathbb Z}[H]} B_{n}H  \]
for $n \geq 0$ by the formula
\[   T_{n}(   1 \otimes_{G} (1, g_{1}, \ldots , g_{n})  ) =   \sum_{i=1}^{m}  \
 1 \otimes_{H} (1, h(i, g_{1}), \ldots , h(i, g_{n}))   .  \]
 Observe that $h(i,1) = 1$.
 \begin{proposition}{$_{}$}
 \label{6.1}
 \begin{em}

 For $n \geq 1$
\[   (1 \otimes d) T_{n} =  T_{n-1}( 1 \otimes d) :  {\mathbb Z} \otimes_{{\mathbb Z}[G]} B_{n}G  \longrightarrow   {\mathbb Z} \otimes_{{\mathbb Z}[H]} B_{n-1}H  .   \]
 \end{em}
 \end{proposition}
 \vspace{2pt}

 {\bf Proof}
 \vspace{2pt}

 We have
 \[ \begin{array}{l}
 (1 \otimes d) T_{n}(  1 \otimes_{G} (1, g_{1}, \ldots , g_{n})   )   \\
 \\
 =   (1 \otimes d)( \sum_{i=1}^{m}  \
 1 \otimes_{H} (1, h(i, g_{1}), \ldots , h(i, g_{n}))  )  \\
 \\
 =    \sum_{i=1}^{m}  \
 1 \otimes_{H} (1, h(i, g_{1})^{-1} h(i, g_{2}) , \ldots ,  h(i, g_{1})^{-1} h(i, g_{n}))   \\
 \\
 \hspace{30pt}  +
  \sum_{i=1}^{m}  \   \sum_{j=1}^{n}  \ (-1)^{j}  \
 1 \otimes_{H} (1, h(i, g_{1}), \ldots , \hat{h(i, g_{j})} \ldots ,  h(i, g_{n})) \\
 \\
 =    \sum_{i=1}^{m}  \
 1 \otimes_{H} (1, h(\pi(g_{1})(i), g_{1}^{-1}g_{2}) , \ldots ,   h(\pi(g_{1})(i), g_{1}^{-1}g_{n}))   \\
 \\
 \hspace{30pt}  +
  \sum_{i=1}^{m}  \   \sum_{j=1}^{n}  \ (-1)^{j}  \
 1 \otimes_{H} (1, h(i, g_{1}), \ldots , \hat{h(i, g_{j})} \ldots ,  h(i, g_{n})) \\
 \\
 =    \sum_{i=1}^{m}  \
 1 \otimes_{H} (1, h(i , g_{1}^{-1}g_{2}) , \ldots ,   h(i , g_{1}^{-1}g_{n}))   \\
 \\
 \hspace{30pt}  +
  \sum_{i=1}^{m}  \   \sum_{j=1}^{n}  \ (-1)^{j}  \
 1 \otimes_{H} (1, h(i, g_{1}), \ldots , \hat{h(i, g_{j})} \ldots ,  h(i, g_{n})) \\
 \\
=  T_{n-1}(1 \otimes d)(  1 \otimes_{G} (1, g_{1}, \ldots , g_{n})   )
 \end{array}  \]
 as required.  $\Box$
 \begin{definition}
 \label{6.2}
 \begin{em}

 The induced homomorphism on homology for $n \geq 0$
 \[   Tr = (T_{n})_{*}  :   H_{n}(G ; {\mathbb Z}) \longrightarrow    H_{n}(H ; {\mathbb Z}) \]
 is called the transfer or corestriction (see \cite{Sn89} p.12).

 This is seen as follows. Define a ${\mathbb Z}[G]$-module chain map
 \[ \tilde{T} :  B_{*}G \longrightarrow  {\mathbb Z}[G]   \otimes_{  {\mathbb Z}[H] }  B_{*}G  \]
 by the formula $\tilde{T}(z) =  \sum_{i=1}^{m} \  x_{i}^{-1} \otimes_{  {\mathbb Z}[H] } x_{i}z$. This chain map, by definition (\cite{Sn89} p.12), induces the transfer
\[ \begin{array}{l}
 \tilde{T} _{*} :   H_{*}(G ; {\mathbb Z}) =  H_{*}( {\mathbb Z}   \otimes_{  {\mathbb Z}[G] }  B_{*}G) \longrightarrow   \\
 \\
 \hspace{40pt}
H_{*}( {\mathbb Z}   \otimes_{  {\mathbb Z}[G] }  {\mathbb Z}[G]   \otimes_{  {\mathbb Z}[H] } B_{*}G)
\cong  H_{*}( {\mathbb Z}    \otimes_{  {\mathbb Z}[H] } B_{*}G) \cong  H_{*}(H ; {\mathbb Z})  .
\end{array}  \]
In this composition the final isomorphism is the inverse to that given by the ${\mathbb Z}[H]$-module chain map
$B_{*}H \longrightarrow  B_{*}G$ induced by the inclusion of $H$ into $G$. However, there is a left inverse to this, namely  the left ${\mathbb Z}[H]$-module chain map
\[    s : B_{*}G \longrightarrow  B_{*}H   \]
given by $s( h_{0}x_{i_{0}} ,  h_{1}x_{i_{1}}  , \ldots , h_{n}x_{i_{n}} ) =  ( h_{0} ,  h_{1}  , \ldots , h_{n})
$ for $h_{j} \in H$. Therefore $s_{*}$ induces on homology the final isomorphism in the composition
\[  H_{*}( {\mathbb Z}    \otimes_{  {\mathbb Z}[H] } B_{*}G) \stackrel{\cong}{\longrightarrow}
H_{*}( {\mathbb Z}    \otimes_{  {\mathbb Z}[H] } B_{*}H)  =   H_{*}(H ; {\mathbb Z})  . \]
It is easy to verify that
\[    T  =   s \cdot \tilde{T} :   {\mathbb Z} \otimes_{{\mathbb Z}[G]} B_{*}G  \longrightarrow   {\mathbb Z} \otimes_{{\mathbb Z}[H]} B_{*}H . \]
 \end{em}
 \end{definition}

\end{document}